\documentclass[a4paper,11pt]{amsart}
\usepackage{graphicx}
 \usepackage{mathptmx}
\usepackage{amsmath}
\usepackage{amssymb}
\usepackage{enumitem}
\usepackage{xcolor}

\newmuskip\pFqmuskip

\newcommand*\pFq[6][8]{%
  \begingroup 
  \pFqmuskip=#1mu\relax
  \mathcode`=\string"8000
  \begingroup\lccode`\~=`\,
  \lowercase{\endgroup\let~}\pFqcomma
  F^{#2}_{#3}{\left(\genfrac..{0pt}{}{#4}{#5}\bigg|#6\right)}%
  \endgroup
}
\newcommand{\pFqcomma}{\mskip\pFqmuskip}

\newtheorem{theorem}{Theorem}
\newtheorem{lemma}{Lemma}

\begin{document}

\title{Some numbers and polynomials related to degenerate harmonic and degenerate hyperharmonic numbers}

\author{Dae San  Kim }
\address{Department of Mathematics, Sogang University, Seoul 121-742, Republic of Korea}
\email{dskim@sogang.ac.kr}

\author{Taekyun  Kim*}
\address{Department of Mathematics, Kwangwoon University, Seoul 139-701, Republic of Korea}
\email{tkkim@kw.ac.kr}

\thanks{ * is corresponding author}

\subjclass[2010]{11B73; 11B83}
\keywords{degenerate harmonic-Fubini polynomials; degenerate hyperharmonic-Fubini polynomials}

\begin{abstract}
Recently, the degenerate harmonic and the degenerate hyperharmonic numbers are introduced respectively as degenerate versions of the harmonic and the hyperharmonic numbers. The aim of this paper is to introduce the degenerate harmonic-Fubini polynomials and numbers related to the degenerate harmonic numbers and to study their properties, explicit expressions and some identities. In addition, as generalizations of those polynomials and numbers, we also introduce the degenerate hyperharmonic-Fubini polynomials and numbers related to the degenerate hyperharmonic numbers and derive similar results to the degenerate harmonic-Fubini polynomials and numbers.
\end{abstract}

\maketitle

\markboth{\centerline{\scriptsize  Some numbers and polynomials related to degenerate harmonic and degenerate hyperharmonic numbers}}
{\centerline{\scriptsize  D. S. Kim and T. Kim}}

\section{Introduction}
It is remarkable that various degenerate versions of many special polynomials and numbers have been studied in recent years with regained interest on them. This exploration for degenerate versions was initiated by Carlitz's work on the degenerate Bernoulli and the degenerate Euler numbers (see [5]). These investigations have been carried out by using such diverse tools as combinatorial methods, generating functions, umbral calculus, $p$-adic analysis, differential equations, probability theory, operator theory, analytic number theory, and so on. \par
The aim of this paper is to introduce the degenerate harmonic-Fubini polynomials (see \eqref{20}) and numbers related to the degenerate harmonic numbers (see \eqref{15}) and to study some properties, explicit expressions and identities for them. As generalizations of those polynomials and numbers, the degenerate hyperharmonic-Fubini polynomials (see \eqref{27}) and numbers related to the degenerate hyperharmonic numbers (see \eqref{16}, \eqref{18}) are also investigated and similar results to the degenerate harmonic-Fubini polynomials and numbers are obtained. \par
The outline of this paper is as follows. In Section 1, we recall the degenerate logarithms together with their properties and the degenerate exponentials. We remind the reader of the degenerate Stirling numbers of the first kind and those of the second kind. We recall the degenerate Fubini polynomials and the generalized degenerate Fubini polynomials. We remind the reader of the degenerate harmonic numbers and the degenerate hyperharmonic numbers. Then we state a useful lemma giving functional equations for two power series. Section 2 is the main result of this paper. We introduce the degenerate harmonic-Fubini polynomials and numbers related to the degenerate harmonic numbers. In Theorem 1, we express the degenerate harmonic-Fubini polynomial a finite sum involving the degenerate Stirling numbers of the second kind and the degenerate harmonic numbers. In Theorem 2, we find an expression of the degenerate harmonic number as a finite sum involving the degenerate harmonic-Fubini numbers and the degenerate Stirling numbers of the first kind. Some generalized degenerate Fubini polynomial is represented in terms of the degenerate Stirling numbers of the second kind in Theorem 3. The degenerate harmonic-Fubini polynomial is expressed as an infinite sum involving the degenerare harmonic numbers in Theorem 4. In Theorem 5, the degenerate hyperharmonic-Fubini polynomial is expressed in terms of the degenerate hyperharmonic numbers and the degenerate Stirling numbers of the second kind, and also in terms of the degenerate harmonic numbers and the degenerate Stirling numbers of the second kind. In Theorem 6, we express the degenerate hyperharmonic-Fubini polynomial as an infinite sum involving the degenerate hyperharmonic numbers. Explicit expressions for the degenerate harmonic numbers are obtained in Theorem 7. In Theorem 8, a functional equation is obtained for any power series $f(t)$ by applying Lemma 1 with $g(t)=-\frac{1}{1-t}\log_{\lambda}(1-t)$. By applying this functional equation to $f(x)=x^{k},\,(k \ge 1)$, we get an identity involving the degenerate harmonic numbers, the degenerate harmonic-Fubini polynomials and some generalized degenerate Fubini polynomials in Theorem 9. From Theorem 9,  an identity of similar nature is derived in Theorem 10. An expression involving the degenerate harmonic-Fubini polynomials and some degenerate Fubini polynomials is shown to be equal to a differential operator applied to $g(x)$, for the aforementioned $g(t)$. Finally, explicit expressions for the degenerate hyperharmonic numbers are founded in Theorem 12. For the rest of this section, we recall the facts that are needed throughout this paper.\par
For any nonzero $\lambda \in\mathbb{R}$, the degenerate logarithms are defined by
\begin{equation}
\log_{\lambda}(1+t)=\sum_{k=1}^{\infty}\frac{(1)_{k,1/\lambda}\lambda^{k-1}}{k!}t^{k}=\frac{1}{\lambda}\big((1+t)^{\lambda}-1\big),\quad (\mathrm{see}\ [4]), \label{1}	
\end{equation}
where
\begin{equation}
(x)_{0,\lambda}=1,\quad (x)_{n,\lambda}=x(x-\lambda)(x-2\lambda)\cdots(x-(n-1)\lambda),\quad (n\ge 1).\label{2}
\end{equation}
From \eqref{1}, we note that
\begin{equation}
\log_{\lambda}(AB)=A^{\lambda}\log_{\lambda}B+\log_{\lambda}A=B^{\lambda}\log_{\lambda}A+\log_{\lambda}B,\label{3}	
\end{equation}
and
\begin{equation}
\log_{\lambda}\bigg(\frac{B}{A}\bigg)=\frac{1}{A^{\lambda}}\Big(\log_{\lambda}B-\log_{\lambda}A\Big),\quad (\mathrm{see}\ [9-13,15]).\nonumber
\end{equation}
For any nonzero $\lambda \in\mathbb{R}$, the degenerate exponentials $e_{\lambda}^{x}(t)$ are defined by
\begin{equation}
e_{\lambda}^{x}(t)=(1+\lambda t)^{\frac{x}{\lambda}}=\sum_{n=0}^{\infty}(x)_{n,\lambda}\frac{t^{n}}{n!},\quad e_{\lambda}(t)=e_{\lambda}^{1}(t),\quad (\mathrm{see}\ [10,14]). \label{4}	
\end{equation}
Note that $\displaystyle\lim_{\lambda\rightarrow 0}\log_{\lambda}(1+t)=\log(1+t),\ \lim_{\lambda\rightarrow 0}e_{\lambda}^{x}(t)=e^{xt}\displaystyle$, and $e_{\lambda}(\log_{\lambda} (t))=\log_{\lambda}(e_{\lambda}(t))=t$.
Thus the  inverse of the degenerate logarithm $\log_{\lambda}(t)$ is the degenerate exponential $e_{\lambda}(t)$. \par
The degenerate Stirling numbers of the first kind are defined  by
\begin{equation}
(x)_{n}=\sum_{k=0}^{n}S_{1,\lambda}(n,k)(x)_{k,\lambda},\quad (n\ge 0),\quad (\mathrm{see}\ [9]), \label{5}
\end{equation}
where $(x)_{0}=1,\ (x)_{n}=x(x-1)\cdots(x-n+1),\ (n\ge 1)$. \par
In addition, the degenerate unsigned Stirling numbers of the first kind are given by
\begin{equation}
{n \brack k}=(-1)^{n-k}S_{1,\lambda}(n,k),\quad (n,k\ge 0),\quad (\mathrm{see}\ [11]). \label{6}
\end{equation}
As the inversion formula of \eqref{5}, the degenerate Stirling numbers of the second kind are given by
\begin{equation}
(x)_{n,\lambda}=\sum_{k=0}^{n}{n \brace k}_{\lambda}(x)_{k},\quad (n\ge 0),\quad (\mathrm{see}\ [9]).\label{7}
\end{equation}
In [15], the degenerate Fubini polynomials are given by
\begin{equation}
\frac{1}{1-x(e_{\lambda}(t)-1)}=\sum_{n=0}^{\infty}F_{n,\lambda}(x)\frac{t^{n}}{n!}.\label{8}
\end{equation}
In particular, for $x=1$, $F_{n,\lambda}=F_{n,\lambda}(1)$ are called the degenerate Fubini numbers. \par
From \eqref{8}, we have
\begin{equation}
F_{n,\lambda}(x)=\sum_{k=0}^{n}{n \brace k}_{\lambda}k!x^{k},\quad (n\ge 0),\quad (\mathrm{see}\ [15]). \label{9}
\end{equation}
Note that $F_{n}(x)=\lim_{\lambda\rightarrow 0}F_{n,\lambda}(x)$ are the ordinary Fubini polynomials given by
\begin{equation}
\frac{1}{1-x(e^{t}-1)}=\sum_{n=0}^{\infty}F_{n}(x)\frac{t^{n}}{n!},\quad (\mathrm{see}\ [1-8,16-20]). \label{10}	
\end{equation}
For any $\alpha\in\mathbb{R}$, the generalized degenerate Fubini polynomials (called the degenerate Fubini polynomials of order $\alpha$) are given by
\begin{equation}
\bigg(\frac{1}{1-x(e_{\lambda}(t)-1)}\bigg)^{\alpha}=\sum_{k=0}^{\infty}F_{k,\lambda}^{(\alpha)}(x)\frac{t^{k}}{k!},\quad (\mathrm{see}\ [15]), \label{11}
\end{equation}
Thus, by \eqref{11}, we get
\begin{equation}
F_{k,\lambda}^{(\alpha)}(x)=\sum_{k=0}^{n}\langle \alpha\rangle_{k}x^{k}{n \brace k}_{\lambda},\quad (\mathrm{see}\ [15]),	\label{12}
\end{equation}
where $\langle \alpha\rangle_{0},\ \langle \alpha\rangle_{k}=\alpha(\alpha+1)\cdots(\alpha+k-1),\ (k\ge 1)$. \par
We recall that the Stirling numbers of the first kind $S_{1}(n,k)$ and those of the second kind ${n \brace k}$ are defined by
\begin{displaymath}
	(x)_{n}=\sum_{k=0}^{n}S_{1}(n,k)x^{k},\quad x^{n}=\sum_{k=0}^{n}{n \brace k}(x)_{k},\quad (n\ge 0),\quad (\mathrm{see}\ [3,4,6,8,12,17]).
\end{displaymath}
Note that $\displaystyle \lim_{\lambda\rightarrow 0}S_{1,\lambda}(n,k)=S_{1}(n,k),\ \lim_{\lambda\rightarrow 0}{n \brace k}_{\lambda}={n\brace k}$. \par
It is well known that the harmonic numbers are defined by
\begin{equation}
H_{0}=0,\quad H_{n}=1+\frac{1}{2}+\cdots+\frac{1}{n},\quad (n\ge 1),\quad (\mathrm{see}\ [6,7,17]). \label{13}
\end{equation}
From \eqref{13}, we note that the generating function of the harmonic numbers is given by
\begin{equation}
-\frac{\log(1-t)}{1-t}=\sum_{k=1}^{\infty}H_{k}t^{k},\quad (\mathrm{see}\ [6,7,17]).\label{14}
\end{equation}
Recently, the degenerate harmonic numbers are defined by
\begin{equation}
-\frac{\log_{\lambda}(1-t)}{1-t}=\sum_{n=1}^{\infty}H_{n,\lambda}t^{n},\quad (\mathrm{see}\ [12,13]). \label{15}
\end{equation}
For $n\ge 0$, $r\ge 1$, the degenerate hyperharmonic numbers are defined by
\begin{equation}
\ H_{0,\lambda}^{(r)}=0\ (r\ge 1),\quad H_{n,\lambda}^{(1)}=H_{n,\lambda}, \quad H_{n,\lambda}^{(r)}=\sum_{k=1}^{n}H_{k,\lambda}^{(r-1)},\ (r\ge 2,\ n\ge 1),\ (\mathrm{see}\ [11,13]). \label{16}
\end{equation}
From \eqref{16}, we note that
\begin{equation}
	H_{n,\lambda}^{(r+1)}=\frac{\binom{n+r}{r}}{\binom{r-\lambda}{r}}(H_{n+r,\lambda}-H_{r,\lambda}),\quad (n,r\in\mathbb{N}), \quad (\mathrm{see}\ [11,13]).\label{17}
\end{equation}
The generating function of the degenerate hyperharmonic numbers is given by
\begin{equation}
-\frac{\log_{\lambda}(1-t)}{(1-t)^{r}}=\sum_{n=1}^{\infty}H_{n,\lambda}^{(r)}t^{n},\quad (r\in\mathbb{N}),\quad (\mathrm{see}\ [11,13]).\label{18}
\end{equation} \par
For $f(x)=\sum_{n=0}^{\infty}a_{n}x^{n}\in\mathbb{C}[\![x]\!]$, we define
\begin{equation*}
f_{\lambda}(x)=\sum_{n=0}^{\infty}a_{n}(x)_{n,\lambda}\in \mathbb{C}[\![x]\!],
\end{equation*}
where $\lambda$ is any fixed real number.
Now, we introduce the next lemma which contains functional equations useful for deriving identities in this paper.
\begin{lemma}[\cite{11}, Theorem 2]
Let $f(x)=\sum_{n=0}^{\infty}a_{n}x^{n},\,\,g(x)=\sum_{k=0}^{\infty}b_{k}x^{k}\in\mathbb{C}[\![x]\!]$. Then we have
\begin{displaymath}
\sum_{n=0}^{\infty}\frac{f^{(n)}(0)}{n!}\sum_{k=0}^{n}{n+r \brace k+r}_{r,\lambda}x^{k}g^{(k)}(x)=\sum_{n=0}^{\infty}\frac{g^{(n)}(0)}{n!}f_{\lambda}(n+r)x^{n}.
\end{displaymath}
In particular, we also have
\begin{equation}
\sum_{m=r}^{\infty}\frac{f^{(m)}(0)}{m!}\bigg(\sum_{k=r}^{m}{m \brace k}_{r,\lambda}x^{k}g^{(k)}(x)\bigg)=\sum_{n=r}^{\infty}\frac{g^{(n)}(0)}{n!}\binom{n}{r}r!\bigg(\sum_{m=r}^{\infty}\frac{f^{(m)}(0)}{m!}(n)_{m-r,\lambda}\bigg)x^{n},\label{19}
\end{equation}
where ${n+r \brace k+r}_{r,\lambda}$ are the degenerate $r$-Stirling numbers defined by
\begin{displaymath}
(x+r)_{n,\lambda}=\sum_{k=0}^{n}{n+r \brace k+r}_{r,\lambda}(x)_{k},\quad (n\ge 0).
\end{displaymath}
\end{lemma}

\section{Some numbers and polynomials related to degenerate harmonic and degenerate hyperharmonic numbers}
In view of \eqref{8} and \eqref{15}, we consider the {\it{degenerate harmonic-Fubini polynomials}} given by
\begin{equation}
-\frac{\log_{\lambda}\big(1-x(e_{\lambda}(t)-1)\big)}{1-x(e_{\lambda}(t)-1)}=\sum_{n=1}^{\infty}HF_{n,\lambda}(x)\frac{t^{n}}{n!}. \label{20}	
\end{equation}
When $x=1$, $HF_{n,\lambda}=HF_{n,\lambda}(1)$ are called the {\it{degenerate harmonic-Fubini numbers}}. \par
From \eqref{14}, we note that
\begin{align}
-\frac{\log_{\lambda}\big(1-x(e_{\lambda}(t)-1)\big)}{1-x(e_{\lambda}(t)-1)}&=\sum_{k=1}^{\infty}H_{k,\lambda}x^{k}(e_{\lambda}(t)-1)^{k} \label{21} \\
&=\sum_{k=1}^{\infty}H_{k,\lambda}x^{k}k!\frac{1}{k!}\big(e_{\lambda}(t)-1\big)^{k} \nonumber \\
&=\sum_{k=1}^{\infty}H_{k,\lambda}k!x^{k}\sum_{n=k}^{\infty}{n \brace k}_{\lambda}\frac{t^{n}}{n!}\nonumber 	\\
&=\sum_{n=1}^{\infty}\sum_{k=1}^{n}{n\brace k}_{\lambda}H_{k,\lambda}k!x^{k}\frac{t^{n}}{n!}. \nonumber
\end{align}
Therefore, by \eqref{20} and \eqref{21}, we obtain the following theorem.
\begin{theorem}
For $n\in\mathbb{N}$, we have
\begin{displaymath}
HF_{n,\lambda}(x)=\sum_{k=1}^{n}{n \brace k}_{\lambda}H_{k,\lambda}k!x^{k}.
\end{displaymath}
In particular, for $x=1$, we get 	\begin{displaymath}
HF_{n,\lambda}=\sum_{k=1}^{n}{n \brace k}_{\lambda}H_{k,\lambda}k!.
\end{displaymath}
\end{theorem}
Replacing $t$ by $\log_{\lambda}(1-t)$ in \eqref{20} and \eqref{21}, we get
\begin{align}
\sum_{n=1}^{\infty}x^{n}H_{n,\lambda}(-1)^{n}t^{n}&=\sum_{k=1}^{\infty}HF_{k,\lambda}(x)\frac{1}{k!}\Big(\log_{\lambda}(1-t)\Big)^{k}\label{22} \\
&=\sum_{k=1}^{\infty}HF_{k,\lambda}(x)\sum_{n=k}^{\infty}	S_{1,\lambda}(n,k)(-1)^{n}\frac{t^{n}}{n!} \nonumber\\
&=\sum_{n=1}^{\infty}(-1)^{n}\sum_{k=1}^{n}HF_{k,\lambda}(x)S_{1,\lambda}(n,k)\frac{t^{n}}{n!}.\nonumber
\end{align}
Therefore, by comparing the coefficients on both sides of \eqref{22}, we obtain the following theorem.
\begin{theorem}
For $n\in\mathbb{N}$, we have
\begin{displaymath}
H_{n,\lambda}x^{n}=\sum_{k=1}^{n}S_{1,\lambda}(n,k)HF_{k,\lambda}(x).
\end{displaymath}
In particular, for $x=1$, we have
\begin{displaymath}
H_{n,\lambda}=\sum_{k=1}^{n}S_{1,\lambda}(n,k)HF_{k,\lambda}.
\end{displaymath}
\end{theorem}
From \eqref{11}, we note that
\begin{align}
\sum_{n=0}^{\infty}F_{n,\lambda}^{(1-\lambda)}(x)\frac{t^{n}}{n!}&=\bigg(\frac{1}{1-x(e_{\lambda}(t)-1)}\bigg)^{1-\lambda}=\sum_{k=0}^{\infty}\binom{k-\lambda}{k}x^{k}(e_{\lambda}(t)-1)^{k} \label{23} \\
&=\sum_{n=0}^{\infty}\langle 1-\lambda\rangle_{k}x^{k}\frac{1}{k!}(e_{\lambda}(t)-1)^{k} \nonumber \\
&=\sum_{k=0}^{\infty}\langle 1-\lambda\rangle_{k}x^{k}\sum_{n=k}^{\infty}{n \brace k}_{\lambda}\frac{t^{n}}{n!}\nonumber \\
&=\sum_{n=0}^{\infty}\sum_{k=0}^{n}\langle 1-\lambda \rangle_{k}x^{k}{n \brace \lambda}_{\lambda}\frac{t^{n}}{n!}. \nonumber
\end{align}
Therefore, by comparing the coefficients on both sides of \eqref{23}, we obtain the following theorem.
\begin{theorem}
Let $n$ be a nonnegative integer. Then we have
\begin{equation}
F_{n}^{(1-\lambda)}(x)=\sum_{k=0}^{n}\langle 1-\lambda\rangle_{k}{n \brace k}_{\lambda}x^{k}.\label{24}
\end{equation}
\end{theorem}
Now, by using \eqref{3} we observe that
\begin{align}
&-\Big(\frac{1}{1-y}\Big) \frac{\log_{\lambda}\Big(1-\frac{y}{1-y}(e_{\lambda}(t)-1)\Big)}{1-\frac{y}{1-y}(e_{\lambda}(t)-1)}=-\frac{\log_{\lambda}\Big(\frac{1-ye_{\lambda}(t)}{1-y}\Big)}{1-ye_{\lambda}(t)}\label{25} \\
&=\frac{-\Big(\frac{1}{1-y}\Big)^{\lambda}\Big(\log_{\lambda}\big(1-ye_{\lambda}(t)\big)-\log_{\lambda}(1-y)\Big)}{1-ye_{\lambda}(t)} \nonumber\\
&=-\Big(\frac{1}{1-y}\Big)^{\lambda}\frac{\log_{\lambda}(1-ye_{\lambda}(t))}{1-ye_{\lambda}(t)}+\frac{\log_{\lambda}(1-y)}{(1-y)^{\lambda}}\frac{1}{1-ye_{\lambda}(t)} \nonumber \\
&=\frac{1}{(1-y)^{\lambda}}\sum_{k=0}^{\infty}H_{k,\lambda}y^{k}e_{\lambda}^{k}(t)+\frac{\log_{\lambda}(1-y)}{(1-y)^{\lambda}}\sum_{k=0}^{\infty}y^{k}e_{\lambda}^{k}(t) \nonumber \\
&=\frac{1}{(1-y)^{\lambda}}\sum_{n=0}^{\infty}\sum_{k=0}^{\infty}\Big(H_{k,\lambda}+\log_{\lambda}(1-y)\Big)y^{k}(k)_{n,\lambda}\frac{t^{n}}{n!}. \nonumber
\end{align}
On the other hand, by \eqref{20}, we get
\begin{equation}
-\Big(\frac{1}{1-y}\Big)\frac{\log_{\lambda}\Big(1-\frac{y}{1-y}(e_{\lambda}(t)-1)\Big)}{1-\frac{y}{1-y}(e_{\lambda}(t)-1)}=\frac{1}{1-y}\sum_{n=0}^{\infty}HF_{n,\lambda}\Big(\frac{y}{1-y}\Big)\frac{t^{n}}{n!}. \label{26}	
\end{equation}
Therefore, by \eqref{25} and \eqref{26}, we obtain the following theorem.

\begin{theorem}
For $n\ge 0$, we have
\begin{displaymath}
\frac{1}{1-y}HF_{n,\lambda}\Big(\frac{y}{1-y}\Big)=\frac{1}{(1-y)^{\lambda}}\sum_{k=0}^{\infty}y^{k}(k)_{n,\lambda}\big(H_{k,\lambda}+\log_{\lambda}(1-y)\big).
\end{displaymath}
Equivalently, we also have
\begin{displaymath}
HF_{n,\lambda}(y)=(1+y)^{\lambda -1}\sum_{k=0}^{\infty}\Big(\frac{y}{1+y}\Big)^{k}(k)_{n,\lambda}\Big(H_{k,\lambda}+\log_{\lambda}\Big(\frac{1}{1+y}\Big)\Big).
\end{displaymath}
Note from Theorem 4 and \eqref{15} that
\begin{displaymath}
\frac{1}{(1-y)^{\lambda}}\bigg\{\bigg(y\frac{d}{dy}\bigg)_{n,\lambda}\Big(-\frac{\log_{\lambda}(1-y)}{1-y}\Big)+\log_{\lambda}(1-y)\bigg(y\frac{d}{dy}\bigg)_{n,\lambda}\frac{1}{1-y}\bigg\}=\frac{1}{1-y}HF_{n,\lambda}\Big(\frac{y}{1-y}\Big).
\end{displaymath}
\end{theorem}
For $r\in\mathbb{N}$, we define the {\it{degenerate hyperharmonic-Fubini polynomials}} given by

\begin{equation}
-\frac{\log_{\lambda}\big(1-y(e_{\lambda}(t)-1\big)}{\big(1-y(e_{\lambda}(t)-1)\big)^{r}}=\sum_{n=1}^{\infty}HF_{n,\lambda}^{(r)}(y)\frac{t^{n}}{n!}. \label{27}
\end{equation}
When $y=1$, $HF_{n,\lambda}^{(r)}=HF_{n,\lambda}^{(r)}(1)$ are called the {\it{degenerate hyperharmonic-Fubini numbers}}. \par
From \eqref{27} and \eqref{17}, we have
\begin{align}
\sum_{n=1}^{\infty}HF_{n,\lambda}^{(r)}(y)\frac{t^{n}}{n!}&=-\frac{\log_{\lambda}(1-y\big(e_{\lambda}(t)-1\big)}{\big(1-y(e_{\lambda}(t)-1)\big)^{r}}	=\sum_{k=1}^{\infty}H_{k,\lambda}^{(r)}\big(e_{\lambda}(t)-1\big)^{k}y^{k} \label{28} \\
&=\sum_{k=1}^{\infty}H_{k,\lambda}^{(r)}y^{k}k!\frac{1}{k!}\big(e_{\lambda}(t)-1\big)^{k}=\sum_{k=1}^{\infty}H_{k,\lambda}^{(r)}y^{k}k!\sum_{n=k}^{\infty}{n \brace k}_{\lambda}\frac{t^{n}}{n!} \nonumber \\
&=\sum_{n=1}^{\infty}\sum_{k=1}^{n}H_{k,\lambda}^{(r)}y^{k}k!{n \brace k}_{\lambda}\frac{t^{n}}{n!} \nonumber \\
&=\sum_{n=1}^{\infty}\Bigg(\sum_{k=1}^{n}\frac{\binom{k+r-1}{r-1}}{\binom{r-1-\lambda}{r-1}}(H_{k+r-1,\lambda}-H_{r-1,\lambda})y^{k}k!{n \brace k}_{\lambda}\bigg)\frac{t^{n}}{n!}. \nonumber
\end{align}
Therefore, by comparing the coefficients on both sides of \eqref{28}, we obtain the following theorem.
\begin{theorem}
Let $n,r$ be positive integers. Then we have
\begin{displaymath}
HF_{n,\lambda}^{(r)}(y)=\sum_{k=1}^{n}H_{k,\lambda}^{(r)}y^{k}k!{n \brace k}_{\lambda}=\sum_{k=1}^{n}\frac{\binom{k+r-1}{r-1}}{\binom{r-1-\lambda}{r-1}}\big(H_{k+r-1,\lambda}-H_{r-1,\lambda}\big)k!{n \brace k}_{\lambda}y^{k}.
\end{displaymath}
In particular, for $y=1$, we get
\begin{displaymath}
HF_{n,\lambda}^{(r)}=\sum_{k=1}^{n}H_{k,\lambda}^{(r)}k!{n \brace k}_{\lambda}=\sum_{k=1}^{n}\frac{\binom{k+r-1}{r-1}}{\binom{r-1-\lambda}{r-1}}\big(H_{k+r-1,\lambda}-H_{r-1,\lambda}\big)k!{n \brace k}_{\lambda}.
\end{displaymath}
\end{theorem}
By using \eqref{3} and \eqref{18}, we observe that
\begin{align}
&-\frac{1}{(1-y)^{r}}\frac{\log_{\lambda}\big(1-\frac{y}{1-y}(e_{\lambda}(t)-1)\big)}{\big(1-\frac{y}{1-y}(e_{\lambda}(t)-1)\big)^{r}}=-\frac{\log_{\lambda}\big(\frac{1-ye_{\lambda}(t)}{1-y}\big)}{\big(1-ye_{\lambda}(t)\big)^{r}} \label{29} \\
&=-\Big(\frac{1}{1-ye_{\lambda}(t)}\Big)^{r}\frac{1}{(1-y)^{\lambda}}\Big(\log_{\lambda}\big(1-ye_{\lambda}(t)\big)-\log_{\lambda}(1-y)\Big) \nonumber \\
&=\frac{1}{(1-y)^{\lambda}}\bigg(-\Big(\frac{1}{1-ye_{\lambda}(t)}\Big)^{r}\log_{\lambda}(1-ye_{\lambda}(t))+\frac{\log_{\lambda}(1-y)}{(1-ye_{\lambda}(t))^{r}}\bigg) \nonumber \\
&=	\frac{1}{(1-y)^{\lambda}}\bigg(\sum_{k=1}^{\infty}H_{k,\lambda}^{(r)}y^{k}e_{\lambda}^{k}(t)+\log_{\lambda}(1-y)\sum_{k=0}^{\infty}\binom{r+k-1}{k}y^{k}e_{\lambda}^{k}(t)\bigg) \nonumber \\
&=\sum_{n=0}^{\infty} \frac{1}{(1-y)^{\lambda}}\sum_{k=0}^{\infty}y^{k}(k)_{n,\lambda}\bigg(H_{k,\lambda}^{(r)}+\binom{r+k-1}{k}\log_{\lambda}(1-y)\bigg)\frac{t^{n}}{n!}.\nonumber
\end{align}
Therefore, by \eqref{27} and \eqref{29}, we obtain the following theorem.
\begin{theorem}
For $n\ge 0$ and $r\in\mathbb{N}$, we have
\begin{displaymath}
\frac{1}{(1-y)^{r}}HF_{n,\lambda}^{(r)}\Big(\frac{y}{1-y}\Big)=\frac{1}{(1-y)^{\lambda}}\sum_{k=0}^{\infty}(k)_{n,\lambda}\bigg(H_{k,\lambda}^{(r)}+\binom{r+k-1}{k}\log_{\lambda}(1-y)\bigg)y^{k}.
\end{displaymath}
Equivalently, we also have
\begin{displaymath}
HF_{n,\lambda}^{(r)}(y)=(1+y)^{\lambda-r}\sum_{k=0}^{\infty}(k)_{n,\lambda}\bigg(H_{k,\lambda}^{(r)}+\binom{r+k-1}{k} \log_{\lambda}\Big(\frac{1}{1+y}\Big)\bigg)\Big(\frac{y}{1+y}\Big)^{k}.
\end{displaymath}
\end{theorem}
Let
\begin{displaymath}
	g(t)=-\frac{\log_{\lambda}(1-t)}{1-t}=-\frac{1}{\lambda}\big((1-t)^{\lambda-1}-(1-t)^{-1}\big), \quad(\mathrm{see}\ \eqref{1}).
\end{displaymath}
Then, from \eqref{1} we have

\begin{align}
g^{(k)}(t)&=\bigg(\frac{d}{dt}\bigg)^{k}g(t)=-\frac{1}{\lambda}\frac{\langle 1-\lambda\rangle_{k}}{(1-t)^{k+1}}(1-t)^{\lambda}+\frac{1}{\lambda}\frac{k!}{(1-t)^{k+1}} \label{30}\\
&=-\frac{\langle 1-\lambda\rangle_{k}}{(1-t)^{k+1}}\log_{\lambda}(1-t)+\frac{1}{(1-t)^{k+1}}\bigg(\frac{k!-\langle 1-\lambda\rangle_{k}}{\lambda}\bigg). \nonumber
\end{align}
From \eqref{30}, we have

\begin{align}
g^{(k)}(0)=\bigg(\frac{d}{dt}\bigg)^{k}g(t)\bigg|_{t=0}=\frac{k!-\langle 1-\lambda\rangle_{k}}{\lambda},\quad (k\in\mathbb{N}). 	\label{31}
\end{align}
By \eqref{15}, we get
\begin{equation}
g^{(k)}(0)=\bigg(\frac{d}{dt}\bigg)^{k}\bigg(-\frac{\log_{\lambda}(1-t)}{1-t}\bigg)\bigg|_{t=0}=\bigg(\frac{d}{dt}\bigg)^{k}\sum_{k=1}^{\infty}H_{k,\lambda}t^{k}\bigg|_{k=0}=k!H_{k,\lambda}.\label{32}
\end{equation}
Therefore, by \eqref{31} and \eqref{32}, we obtain the following theorem.
\begin{theorem}
For $k\in\mathbb{N}$, we have
\begin{displaymath}
H_{k,\lambda}=\frac{1}{k!}\bigg(\frac{k!-\langle 1-\lambda \rangle_{k}}{\lambda}\bigg)=\frac{1}{\lambda}\bigg(1-\binom{k-\lambda}{k}\bigg).
\end{displaymath}
\end{theorem}
Let $r$ be a nonnegative integer, $f(x) \in\mathbb{C}[\![x]\!]$, and let $g(t)=-\frac{1}{1-t}\log_{\lambda}(1-t)$.\\
By \eqref{19}, \eqref{30} and \eqref{32}, we get
\begin{align}
	&\sum_{n=r}^{\infty}\binom{n}{r}H_{n,\lambda}r!\bigg(\sum_{m=r}^{\infty}\frac{f^{(m)}(0)}{m!}(n)_{m-r,\lambda}\bigg)x^{n}\label{33}\\
	&=\frac{1}{1-x}\sum_{m=r}^{\infty}\frac{f^{(m)}(0)}{m!}\sum_{k=r}^{m}{m \brace k}_{r,\lambda}\bigg(\frac{k!-\langle 1-\lambda\rangle_{k}}{\lambda}\bigg)\Big(\frac{x}{1-x}\Big)^{k} \nonumber \\
	&\quad -\frac{1}{1-x}\log_{\lambda}(1-x)\sum_{m=r}^{\infty}\frac{f^{(m)}(0)}{m!}\sum_{k=r}^{m}{m \brace k}_{r,\lambda}\langle 1-\lambda\rangle_{k}\Big(\frac{x}{1-x}\Big)^{k}. \nonumber
\end{align}
Let $r=0$ in \eqref{33}. Then, by Theorem1, Theorem 7 and \eqref{24}, we have
\begin{align}
&\sum_{n=0}^{\infty}H_{n,\lambda} \bigg(\sum_{m=0}^{\infty}\frac{f^{(m)}(0)}{m!}(n)_{m,\lambda}\bigg)x^{n} \label{34} \\
&=\frac{1}{1-x}\sum_{m=0}^{\infty}\frac{f^{(m)}(0)}{m!}\sum_{k=0}^{m}{m \brace k}_{\lambda}\bigg(\frac{k!-\langle 1-\lambda\rangle_{k}}{\lambda}\bigg)\Big(\frac{x}{1-x}\Big)^{k} \nonumber \\
&\quad -\frac{1}{1-x}\log_{\lambda}(1-x)\sum_{m=0}^{\infty}\frac{f^{(m)}(0)}{m!}\sum_{k=0}^{m}{m \brace k}_{\lambda}\langle 1-\lambda\rangle_{k}\Big(\frac{x}{1-x}\Big)^{k}\nonumber \\
&=\frac{1}{1-x}\sum_{m=0}^{\infty}\frac{f^{(m)}(0)}{m!}\sum_{k=0}^{m}{m\brace k}_{\lambda}k!H_{k,\lambda}\Big(\frac{x}{1-x}\Big)^{k}-\frac{\log_{\lambda}(1-x)}{1-x}\sum_{n=0}^{\infty}\frac{f^{(m)}(0)}{m!}F_{m,\lambda}^{(1-\lambda)}\Big(\frac{x}{1-x}\Big)\nonumber \\
&=\frac{1}{1-x}\sum_{m=0}^{\infty}\frac{f^{(m)}(0)}{m!}HF_{m,\lambda}\Big(\frac{x}{1-x}\Big)-\frac{\log_{\lambda}(1-x)}{1-x}\sum_{m=0}^{\infty}\frac{f^{(m)}(0)}{m!}F_{m,\lambda}^{(1-\lambda)}\Big(\frac{x}{1-x}\Big). 	\nonumber
\end{align}
Therefore, by \eqref{34}, we obtain the following theorem.
\begin{theorem}
Let $f(x) \in\mathbb{C}[\![x]\!]$. Then we have
\begin{equation}
\begin{aligned}
&\sum_{n=0}^{\infty}H_{n,\lambda}\bigg(\sum_{m=0}^{\infty}\frac{f^{(m)}(0)}{m!}(n)_{m,\lambda}\bigg)x^{n}\\
&=\frac{1}{1-x}\sum_{m=0}^{\infty}\frac{f^{(m)}(0)}{m!}HF_{m,\lambda}\Big(\frac{x}{1-x}\Big)-\frac{\log_{\lambda}(1-x)}{1-x}\sum_{m=0}^{\infty}\frac{f^{(m)}(0)}{m!}F_{m,\lambda}^{(1-\lambda)}\Big(\frac{x}{1-x}\Big). \label{35}
\end{aligned}
\end{equation}
\end{theorem}
Let $f(x)=x^{k}, \ (k\ge 1)$,  in \eqref{35}. Then we have
\begin{equation}
\sum_{n=0}^{\infty}H_{n,\lambda}(n)_{k,\lambda}x^{n}=\frac{1}{1-x}HF_{k,\lambda}\Big(\frac{x}{1-x}\Big)-\frac{\log_{\lambda}(1-x)}{1-x}F_{k,\lambda}^{(1-\lambda)}\Big(\frac{x}{1-x}\Big).	\label{36}
\end{equation}

Therefore, by \eqref{35}, we obtain the following theorem.
\begin{theorem}
	For $k\ge 1$, we have
	\begin{displaymath}
		\sum_{n=1}^{\infty}H_{n,\lambda}(n)_{k,\lambda}x^{n}=\frac{1}{1-x}HF_{k,\lambda}\Big(\frac{x}{1-x}\Big)-\frac{\log_{\lambda}(1-x)}{1-x}F_{k,\lambda}^{(1-\lambda)}\Big(\frac{x}{1-x}
\Big).
\end{displaymath}
\end{theorem}
From \eqref{36}, we note that
\begin{align}
&\frac{1}{(1-x)^{2}}\bigg(HF_{k,\lambda}\Big(\frac{x}{1-x}\Big)-\log_{\lambda}(1-x)F_{k,\lambda}^{(1-\lambda)}\Big(\frac{x}{1-x}\Big)\bigg)=\frac{1}{1-x}\sum_{l=1}^{\infty}H_{l,\lambda}(l)_{k,\lambda}x^{l} \label{37}\\
&=\sum_{j=0}^{\infty}x^{j}\sum_{l=1}^{\infty}H_{l,\lambda}(l)_{k,\lambda}x^{l}=\sum_{n=1}^{\infty}\bigg(\sum_{l=1}^{n}H_{l,\lambda}(l)_{k,\lambda}\bigg)x^{n}.\nonumber
\end{align}
Therefore, by \eqref{37}, we obtain the following theorem.
\begin{theorem}
For $k\in\mathbb{N}$, we have
\begin{align*}
&\sum_{n=1}^{\infty}\bigg((1)_{k,\lambda}H_{1,\lambda}+(2)_{k,\lambda}H_{2,\lambda}+\cdots+(n)_{k,\lambda}H_{n,\lambda}\bigg)x^{n}\\
&\quad =\frac{1}{(1-x)^{2}}\bigg(HF_{k,\lambda}\Big(\frac{x}{1-x}\Big)-\log_{\lambda}(1-x)F_{k,\lambda}^{(1-\lambda)}\Big(\frac{x}{1-x}\Big)\bigg).
\end{align*}
\end{theorem}
By \eqref{36}, we get
\begin{equation}
\begin{aligned}
&\sum_{n=1}^{\infty}\big((n)_{1,\lambda}+(n)_{2,\lambda}+\cdots+(n)_{k,\lambda}\big)H_{n,\lambda}x^{n}\\
&=\frac{1}{1-x}\sum_{l=1}^{k}\bigg(HF_{l,\lambda}\Big(\frac{x}{1-x}\Big)-\log_{\lambda}(1-x)F_{l,\lambda}^{(1-\lambda)}\Big(\frac{x}{1-x}\Big)\bigg).
\end{aligned}\label{38}
\end{equation}
From \eqref{15} and \eqref{36}, we note that
\begin{align}
\Big(x \frac{d}{dx}\Big)_{k,\lambda}\bigg(-\frac{\log_{\lambda}(1-x)}{1-x}\bigg)&=\Big(x \frac{d}{dx}\Big)_{k,\lambda}\bigg(\sum_{n=1}^{\infty}H_{n,\lambda}x^{n}\bigg)=\sum_{n=1}^{\infty}(n)_{k,\lambda}H_{n,\lambda}x^{n} \label{39} \\
&=\frac{1}{1-x}\bigg(HF_{k,\lambda}\Big(\frac{x}{1-x}\Big)-\log_{\lambda}(1-x)F_{k,\lambda}^{(1-\lambda)}\Big(\frac{x}{1-x}\Big)\bigg).\nonumber
\end{align} \par
Therefore, by \eqref{39}, we obtain the following differential equation.
\begin{theorem}
Let $k$ be a positive integer. Then we have
\begin{displaymath}
\Big(x \frac{d}{dx}\Big)_{k,\lambda}\bigg(-\frac{\log_{\lambda}(1-x)}{1-x}\bigg)=\frac{1}{1-x}\bigg(HF_{k,\lambda}\Big(\frac{x}{1-x}\Big)-\log_{\lambda}(1-x)F_{k,\lambda}^{(1-\lambda)}\Big(\frac{x}{1-x}\Big)\bigg).
\end{displaymath}
\end{theorem}
For an integer $r$ with $r >1$, we let
\begin{equation}
g(t)=-\frac{1}{(1-t)^{r}}\log_{\lambda}(1-t)=-\frac{1}{\lambda}\big((1-t)^{\lambda-r}-(1-t)^{-r}\big),\quad (\mathrm{see} \,\, \eqref{1}).\label{41}
\end{equation}
Then, for $k\in\mathbb{N}$, we have
\begin{align}
g^{(k)}(t)&=\bigg(\frac{d}{dt}\bigg)^{k}g(t)=-\frac{\langle r-\lambda\rangle_{k}}{\lambda}(1-t)^{\lambda-r-k}+\frac{1}{\lambda}\langle r\rangle_{k}(1-t)^{-r-k} \label{42}\\
&=-\frac{\langle r-\lambda\rangle_{k}}{(1-t)^{r+k}}\log_{\lambda}(1-t)+\frac{1}{(1-t)^{r+k}}\bigg(\frac{\langle r\rangle_{k}-\langle r-\lambda \rangle_{k}}{\lambda}\bigg). \nonumber	
\end{align}
From \eqref{18}, we note that
\begin{equation}
g^{(k)}(0)=k!H_{k,\lambda}^{(r)},\quad (k\ge 1). \label{43}
\end{equation}
Therefore, by \eqref{42} and \eqref{43}, we obtain the following theorem.
\begin{theorem}
For $k\ge 1$, we have
\begin{displaymath}
H_{k,\lambda}^{(r)}=\frac{1}{k!}\bigg(\frac{\langle r\rangle_{k}-\langle r-\lambda\rangle_{k}}{\lambda}\bigg)=\frac{1}{\lambda}\bigg(\binom{r+k-1}{k}-\binom{r+k-\lambda-1}{k}\bigg).
\end{displaymath}
\end{theorem}

Let $f(x)=\sum_{n=0}^{\infty}a_{n}x^{n}\in\mathbb{C}[\![x]\!]$, and let $g(t)=-\frac{\log_{\lambda}(1-t)}{(1-t)^{r}}$. Then, by Theorem 12, \eqref{19}, \eqref{41} and \eqref{42}, we get
\begin{align}
	&\sum_{n=r}^{\infty}H_{n,\lambda}^{(r)}\binom{n}{r}r!\bigg(\sum_{m=r}^{\infty}\frac{f^{(m)}(0)}{m!}(n)_{m-r,\lambda}\bigg)x^{n} \label{45} \\
	&=\sum_{m=r}^{\infty}\frac{f^{(m)}(0)}{m!}\sum_{k=r}^{m}{m\brace k}_{r,\lambda}x^{k}\bigg(-\frac{\langle r-k\rangle_{k}}{(1-x)^{r+k}}\log_{\lambda}(1-x)+\frac{1}{(1-x)^{r+k}}\bigg(\frac{\langle r\rangle_{k}-\langle r-\lambda\rangle_{k}}{\lambda}\bigg)\bigg) \nonumber\\
	&=-\frac{\log_{\lambda}(1-x)}{(1-x)^{r}}\sum_{m=r}^{\infty}\frac{f^{(m)}(0)}{m!}\sum_{k=r}^{m}{m \brace k}_{r,\lambda}\langle r-\lambda\rangle_{k}\Big(\frac{x}{1-x}\Big)^{k}\nonumber \\
	&\qquad +\frac{1}{(1-x)^{r}}\sum_{m=r}^{\infty}\frac{f^{(m)}(0)}{m!}\sum_{k=r}^{m}{m\brace k}_{r,\lambda}\Big(\frac{x}{1-x}\Big)^{k}k!H_{k,\lambda}^{(r)}\nonumber
\end{align}
Let $f(x)=x^{k},\ (k\ge 1)$ in \eqref{45}. Then we have
\begin{equation}
\begin{aligned}
	&\sum_{n=r}^{\infty}H_{n,\lambda}^{(r)}\binom{n}{r}r!(n)_{k-r,\lambda}x^{n}\\
	&=-\frac{\log_{\lambda}(1-x)}{(1-x)^{r}}\sum_{l=r}^{k}{k \brace l}_{r,\lambda}\langle r-  \lambda \rangle _{l} \Big(\frac{x}{1-x}\Big)^{l}+\frac{1}{(1-x)^{r}}\sum_{l=r}^{k}{k \brace l}_{r,\lambda}\Big(\frac{x}{1-x}\Big)^{l}l!H_{l,\lambda}^{(r)}.
\end{aligned}\label{46}
\end{equation}

\section{Conclusion}
The degenerate harmonic-Fubini polynomials are given by $HF_{n,\lambda}(x)=\sum_{k=1}^{n}{n \brace k}_{\lambda}H_{k,\lambda}k!x^{k}$, with $H_{k,\lambda}$ the degenerate harmonic numbers, while the degenerate Fubini polynomials are given by $F_{n,\lambda}(x)=\sum_{k=0}^{n}{n \brace k}_{\lambda}k!x^{k}$. The degenerate harmonic-Fubini polynomials are so named for this reason. The degenerate hyperharmonic-Fubini polynomials are also so named, as it is given by $HF_{n,\lambda}^{(r)}(x)=\sum_{k=1}^{n}{n \brace k}_{\lambda}H_{k,\lambda}^{(r)}k!x^{k}$, with $H_{k,\lambda}^{(r)}$ the degenerate hyperharmonic numbers. \par
We introduced the degenerate harmonic-Fubini polynomials and numbers and studied their properties, explicit expressions and some identities by using generating functions. In addition, as generalizations of those polynomials and numbers, we also introduced the degenerate hyperharmonic-Fubini polynomials and derived similar results to the degenerate harmonic-Fubini polynomials and numbers. \par
It is one of our future projects to continue to study various degenerate versions of some special polynomials and numbers and to find their applications to physics, science and engineering as well as to mathematics.

\end{document}